\documentclass[reqno,12pt]{amsart}
\usepackage{amssymb}
\usepackage{amsfonts}
\usepackage{amsmath}
\usepackage{graphicx}
\usepackage{amscd,amsthm}

\newtheorem{theorem}{Theorem}
\theoremstyle{plain}
\newtheorem{acknowledgement}{Acknowledgement}

\newtheorem{lemma}{Lemma}

\numberwithin{equation}{section}

\input{tcilatex}

\begin{document}
\author{}
\title{}
\maketitle

\begin{center}
\thispagestyle{empty} \pagestyle{myheadings} 
\markboth{\bf Burak Kurt and
Yilmaz Simsek}{\bf On The Hermite Based-Second Kind Genocchi polynomials}

\textbf{\Large On The Hermite Based-Second Kind Genocchi Polynomials}

\bigskip

\textbf{Burak Kurt} \textbf{and Yilmaz Simsek}\\[0pt]
\bigskip

\medskip

Department of Mathematics, Faculty of Science University of AkdenizTR-07058
Antalya, Turkey \textbf{\ }burakkurt@akdeniz.edu.tr and
ysimsek@akdeniz.edu.tr\\[0pt]

\medskip

\textbf{The 24th Congress of the Jangjeon Mathematical Society 21 July
Konya-TURKEY}

\textbf{{\large {Abstract}}}\medskip
\end{center}

\begin{quotation}
The aim of this paper is to study generating function of the Hermite-Kamp%
\.{e} de F\.{e}riet based second kind Genocchi polynomials. We also give
some identities related to these polynomials.
\end{quotation}



\section{Introduction, Definitions and Notations}

The Bernoulli numbers the Bernoulli polynomials, the Euler numbers, the
Euler polynomials and Genocchi polynomials and numbers are used many areas
of sciences. These numbers and polynomials are also used in calculus of
finite differences, in numerical analysis, in analytical number theory.

Recently many authors have studied on these numbers and polynomials. They
also defined many different generating functions of these polynomials and
numbers.

Due to Bretti and Ricci \cite{bretti1}, the Hermite-Kamp%
%TCIMACRO{\U{b4}}%
%BeginExpansion
\'{}%
%EndExpansion
e de F%
%TCIMACRO{\U{b4}}%
%BeginExpansion
\'{}%
%EndExpansion
eriet (or Gould-Hopper) polynomials have been used in order to construct
addition formulas for different classes of generalized Gegenbauer
polynomials. These polynomials (two-variable Hermite polynomials) are
defined by means of the following generating functions:%
\begin{equation}
e^{xt+yt^{j}}=\dsum\limits_{n=0}^{\infty }H_{n}^{\left( j\right) }\left(
x,y\right) \frac{t^{n}}{n!}.  \label{11}
\end{equation}

From (\ref{11}), one can easily see that%
\begin{equation*}
H_{n}^{(j)}\left( x,y\right) =n!\sum_{r=0}^{[\frac{n}{j}]}\frac{x^{n-jr}y^{r}%
}{r!\left( n-jr\right) !},
\end{equation*}%
where $j\geq 2$ is an integer cf. \cite{bretti1}. In \cite{bretti1}, the
case $j=1$ is not considered, since the corresponding $2D$ polynomials are
simply expressed by the Newton binomial formula.

The polynomials $H_{n}^{(j)}\left( x,y\right) $ are the solution of the
generalized heat equation:%
\begin{eqnarray*}
\frac{\partial }{\partial x}F(x,y) &=&\frac{\partial ^{j}}{\partial x^{j}}%
F(x,y), \\
F(x,0) &=&x^{n}.
\end{eqnarray*}

The classical Genocchi numbers $G_{n}$ and the classical Genocchi numbers $%
G_{n}(x)$ are usually defined by means of the following generating function,%
\begin{equation}
g(t)=\frac{2t}{e^{t}+1}=\sum_{n=0}^{\infty }G_{n}\frac{t^{n}}{n!}\text{, }%
\left\vert t\right\vert <\pi  \label{1}
\end{equation}%
and 
\begin{equation}
g(t,x)=g(t)e^{xt}=\sum_{n=0}^{\infty }G_{n}(x)\frac{t^{n}}{n!}\text{, }%
\left\vert t\right\vert <\pi \text{.}  \label{2}
\end{equation}

From (\ref{1}) and (\ref{2}), one can easily see that%
\begin{equation*}
G_{n}(x)=\sum_{k=0}^{n}\binom{n}{k}G_{k}x^{n-k}\text{.}
\end{equation*}

The second kind Genocchi polynomials of higher order are defined by means of
the following generating function: 
\begin{equation}
\mathfrak{g}(t,x)=\left( \frac{2t}{e^{t}+e^{-t}}\right)
^{a}e^{tx}=\sum_{n=0}^{\infty }\mathcal{G}_{n}^{(a)}(x)\frac{t^{n}}{n!}\text{%
, }\left\vert t\right\vert <\frac{\pi }{2}.  \label{4}
\end{equation}%
Note that $\mathcal{G}_{n}^{(1)}(x)=\mathcal{G}_{n}(x)$ denotes the second
kind Genocchi polynomials cf. \cite{ryoo13}. $\mathcal{G}_{n}(0)=\mathcal{G}%
_{n}$ which denote so-called the second kind Genocchi numbers cf. \cite%
{ryoo13}.

By using (\ref{4}), one can easily see that%
\begin{equation*}
\mathcal{G}_{n}(x)=\sum_{k=0}^{n}\binom{n}{k}x^{n-k}\mathcal{G}_{k}.
\end{equation*}%
By (\ref{2}) and (\ref{4}), we have%
\begin{eqnarray*}
\mathcal{G}_{n}(x) &=&2^{n-1}G_{n}(\frac{x+1}{2}) \\
&=&\sum_{k=0}^{n}\sum_{j=0}^{n-k}\binom{n}{k,j,n-k-j}2^{k-1}x^{n-k-j}G_{k},
\end{eqnarray*}%
where%
\begin{equation*}
\binom{n}{a,b,c}=\frac{n!}{a!b!c!},a+b+c=n.
\end{equation*}

The second kind Euler polynomials of higher order are defined by means of
the following generating function:%
\begin{equation}
\left( \frac{2}{e^{t}+e^{-t}}\right) ^{a}e^{xt}=\sum_{n=0}^{\infty }\mathcal{%
E}_{n}^{(a)}(x)\frac{t^{n}}{n!},\text{ }\left\vert t\right\vert <\frac{\pi }{%
2}  \label{8}
\end{equation}%
Observe that $\mathcal{E}_{n}^{(1)}(x)=\mathcal{E}_{n}(x)$ denotes the
second kind Euler polynomials cf. (\cite{jang10}, \cite{ryoo14}).

Relation between $\mathcal{E}_{n}(x)$ and $\mathcal{G}_{n}(x)$ is given by%
\begin{equation*}
\mathcal{G}_{n}(x)=n\mathcal{E}_{n-1}(x).
\end{equation*}

\begin{lemma}
The second kind Genocchi polynomials of $\alpha $ order is satisfied the
following relations%
\begin{equation}
G_{n}^{\left( \alpha +\beta \right) }(x+x_{1})=\sum_{k=0}^{n}\binom{n}{k}%
G_{k}^{\left( \alpha \right) }(x)G_{n-k}^{\left( \beta \right) }(x_{1}).
\label{9}
\end{equation}%
Also 
\begin{equation}
\sum_{k=0}^{n}\binom{n}{k}E_{k}^{\left( \alpha \right) }G_{n-k}^{\left(
\alpha \right) }(x+x_{1})=\sum_{k=0}^{n}\binom{n}{k}E_{k}^{\left( \alpha
\right) }\left( x\right) G_{n-k}^{\left( \alpha \right) }(x_{1})\text{.}
\label{10}
\end{equation}
\end{lemma}

\begin{proof}
From (\ref{4}),%
\begin{eqnarray*}
\sum_{n=0}^{\infty }G_{n}^{\left( \alpha +\beta \right) }(x+x_{1})\frac{t^{n}%
}{n!} &=&\left( \sum_{n=0}^{\infty }G_{n}^{\left( \alpha \right) }(x)\frac{%
t^{n}}{n!}\right) \left( \sum_{n=0}^{\infty }G_{n}^{\left( \beta \right)
}(x_{1})\frac{t^{n}}{n!}\right) \\
&=&\sum_{n=0}^{\infty }\left( \sum_{k=0}^{n}\binom{n}{k}G_{k}^{\left( \alpha
\right) }(x)G_{n-k}^{\left( \beta \right) }(x_{1})\right) \frac{t^{n}}{n!}%
\text{.}
\end{eqnarray*}%
Comparing the coefficients, we prove Lemma 1.
\end{proof}

\section{Hermite-Kamp\.{e} de F\.{e}riet based second kind Genocchi
polynomials}

The aim of this section is to define the Hermite-Kamp\.{e} de F\.{e}riet (or
Gould-Hopper) based second kind Genocchi polynomials. Some poroperties of
these polynomials are also given.

The Hermite-Kamp\.{e} de F\.{e}riet based second kind Genocchi polynomials
of higher order are defined by means of the following generating function:%
\begin{equation}
F_{H,a}(t,x,y;j)=\left( \frac{2t}{e^{t}+e^{-t}}\right)
^{a}e^{xt+yt^{j}}=\dsum\limits_{n=0}^{\infty }\left( _{H}\mathcal{G}%
_{n}^{(j,a)}\left( x,y\right) \right) \frac{t^{n}}{n!}.  \label{b1}
\end{equation}%
From (\ref{b1}), one finds that%
\begin{equation*}
_{H}\mathcal{E}_{n-1}^{(j,a)}\left( x,y\right) =\frac{_{H}\mathcal{G}%
_{n}^{(j,a)}\left( x,y\right) }{n},
\end{equation*}%
where $_{H}\mathcal{E}_{n-1}^{(j,a)}\left( x,y\right) $ denotes the
Hermite-Kamp\.{e} de F\.{e}riet based second kind Euler polynomials of
higher order.

Now we compute the derivative of ((\ref{b1}) with respect to $x$ to derive a
derivative formula for the Hermite-Kamp\.{e} de F\.{e}riet based second kind
Genocchi polynomials of higher order:%
\begin{equation*}
\frac{\partial }{\partial x}\left( _{H}\mathcal{G}_{n}^{(j,a)}\left(
x,y\right) \right) =n\left( _{H}\mathcal{G}_{n-1}^{(j,a)}\left( x,y\right)
\right) .
\end{equation*}

If we substitute $j=2$ and $a=1$ into (\ref{b1}), and using (\ref{4}), we
have 
\begin{equation}
\dsum\limits_{n=0}^{\infty }\left( _{H}\mathcal{G}_{n}^{(2,1)}\left(
x,y\right) \right) \frac{t^{n}}{n!}=\sum_{n=0}^{\infty }\mathcal{G}_{n}(x)%
\frac{t^{n}}{n!}\sum_{n=0}^{\infty }\frac{y^{n}t^{2n}}{n!}.  \label{12}
\end{equation}

After some elementary calculations in the above, we obtain%
\begin{equation*}
\dsum\limits_{n=0}^{\infty }\left( _{H}\mathcal{G}_{n}^{(2)}\left(
x,y\right) \right) \frac{t^{n}}{n!}=\dsum\limits_{n=0}^{\infty }\left(
n!\sum_{l=0}^{\left[ \frac{n}{2}\right] }\frac{y^{l}\mathcal{G}_{n-2l}(x)}{%
l!(n-2l)!}\right) \frac{t^{n}}{n!}.
\end{equation*}%
By comparing coefficients of $t^{n}$ in the both sides of the above, we
arrive at the result of the theorem:

\begin{theorem}
\begin{equation*}
_{H}\mathcal{G}_{n}^{(2,1)}\left( x,y\right) =n!\sum_{l=0}^{\left[ \frac{n}{2%
}\right] }\frac{y^{l}\mathcal{G}_{n-2l}(x)}{l!(n-2l)!}.
\end{equation*}
\end{theorem}

\begin{theorem}
\begin{equation}
_{H}\mathcal{G}_{n}^{(j,a+b)}\left( x_{1}+x_{2},y_{1}+y_{2}\right)
=\sum_{k=0}^{n}\binom{n}{k}\left( _{H}\mathcal{G}_{n}^{(j,a)}\left(
x_{1},y_{1}\right) \right) \left( _{H}\mathcal{G}_{n-k}^{(j,b)}\left(
x_{2},y_{2}\right) \right) \text{.}  \label{13}
\end{equation}
\end{theorem}

\begin{proof}
By (\ref{b1}), we define the following functional equation%
\begin{equation*}
F_{H,a+b}(t,x_{1}+x_{2},y_{1}+y_{2};j)=F_{H,a}(t,x_{1},y_{1};j)F_{H,b}(t,x_{2},y_{2};j)
\end{equation*}%
By using the above functional equation, we obtain%
\begin{equation*}
\dsum\limits_{n=0}^{\infty }\left( _{H}\mathcal{G}_{n}^{(j,a+b)}\left(
x_{1}+x_{2},y_{1}+y_{2}\right) \right) \frac{t^{n}}{n!}=\dsum\limits_{n=0}^{%
\infty }\left( _{H}\mathcal{G}_{n}^{(j,a)}\left( x_{1},y_{1}\right) \right) 
\frac{t^{n}}{n!}\dsum\limits_{n=0}^{\infty }\left( _{H}\mathcal{G}%
_{n}^{(j,b)}\left( x_{2},y_{2}\right) \right) \frac{t^{n}}{n!}
\end{equation*}%
By using Cauchy product in the above equation, than comparing the
coefficients of $\frac{t^{n}}{n!}$ on both sides, we arrive at the desired
result.
\end{proof}

Relation between Hermite-based second kind Genocchi polynomials of higher
order, the second kind Genocchi polynomials of higher order and two-variable
generalized Hermite polynomials is given by the next theorem.

\begin{theorem}
\begin{equation*}
_{H}\mathcal{G}_{n}^{(2,a)}(x,y)=\sum_{l=0}^{n}\binom{n}{l}\mathcal{G}%
_{n-l}^{\left( a\right) }H_{l}^{\left( 2\right) }\left( x,y\right) .
\end{equation*}
\end{theorem}

\begin{proof}
If we substitute $j=2$ and $a=1$ into (\ref{b1}), we easily arive at the
desired result.
\end{proof}

\begin{acknowledgement}
The present investigation was supported, by the \textit{Scientific Research
Project Administration of Akdeniz University}.
\end{acknowledgement}

\end{document}